\documentclass[11pt,letterpaper]{article}
\usepackage{amsmath}
\usepackage{amsthm}
\usepackage{amsfonts}
\usepackage{comment}
\usepackage{graphicx}
\usepackage{color}

\newtheorem{lemma}{Lemma}

\newtheorem{theorem}{Theorem}

\newcommand{\E}{\mathrm{E}}
\newcommand{\Def}{\overset{\triangle}=}

\def\bfPhi{\mbox{\protect\boldmath$\Phi$}}
\newcommand{\PR}{\mbox{\rm Pr}}

\newlength{\noteWidth}
\long\def\notes#1{\ifinner
             {\tiny #1}
             \else
             \marginpar{\parbox[t]{\noteWidth}{\raggedright\tiny #1}}
             \fi}

\begin{document}

\begin{titlepage}
\title{\Huge Steady State Analysis of Balanced-Allocation Routing}

\author{{\Large Aris Anagnostopoulos}\thanks{
Computer Science Department, Brown University, Box 1910, Providence,
RI 02912-1910, USA. \hbox{E-mail}:~{\tt \{aris,~eli\}@cs.brown.edu.}
Supported in part by NSF grants CCR-0121154, and DMI-0121495.
}
\and
{\Large Ioannis Kontoyiannis}\thanks{Division of Applied Mathematics
and Department of Computer Science, 
Brown University, 
Box F, 182 George St., Providence, RI 02912, USA.
\hbox{E-mail}:~{\tt yannis@cs.brown.edu}
\hbox{Web}:~{\tt www.dam.brown.edu/people/yannis/.}
Supported in part by NSF grant \#0073378-CCR 
and USDA-IFAFS grant \#00-52100-9615.
}
\and
{\Large Eli Upfal}\addtocounter{footnote}{-2}\footnotemark
}

\maketitle
\begin{abstract}
We compare the long-term, steady-state performance of a
variant of the standard \emph{Dynamic Alternative Routing (DAR)}
technique commonly used in telephone and ATM networks,
to the performance of a path-selection algorithm 
based on the ``balanced-allocation'' principle~\cite{abku-ba-00,m-ptcrlb-96};
we refer to this new algorithm as the \emph{Balanced Dynamic 
Alternative Routing (BDAR)} algorithm.
While DAR checks alternative routes sequentially until 
available bandwidth is found, the BDAR algorithm compares 
and chooses the best among a small number of alternatives.

We show that, at the expense of a minor increase in routing 
overhead, the BDAR algorithm gives a substantial improvement 
in network performance, in terms both of network congestion 
and of bandwidth requirement.
\end{abstract}
\thispagestyle{empty}
\end{titlepage}

\section{Introduction}
\label{sec:int}
Fast, high bandwidth, circuit switching telecommunications 
systems such as ATM and telephone networks often employ 
a limited path-selection algorithm in order to fully 
utilize the network resources while minimizing routing 
overhead. Typically,
between each pair of nodes in the network there is 
a dedicated bandwidth for communication, namely, 
no more than a certain fixed number of calls can
be simultaneously active between each pair of nodes.
This dedicated bandwidth is chosen in order to satisfy 
the demand for communication between these stations. 
Only when this bandwidth is exhausted the admission 
control protocol tries to find an alternative route
through intermediate nodes. To minimize overhead
and routing delays, the protocol checks just a 
small number of alternative routes; if there are 
no free connections available on any of these 
alternatives, then the call or communication 
request is rejected.
Implementations that use this technique include
the Dynamic Alternate Routing (DAR) algorithm used by
British Telecom~\cite{gkk-dar-95},
and AT\&T's Dynamic Nonhierarchical Routing (DNHR)
algorithm~\cite{acm-dondr-81}.

A common feature in these (and other)
currently implemented protocols 
is the sequential examination of alternative routes. 
Only when the algorithm examines a route and finds it
cannot be used an alternative one is examined. 
The criteria for when a route can or should
be used, and the method in which the alternative 
route is selected have been the subject of 
extensive research, in particular, 
in the context of British Telecom's DAR 
algorithm~\cite{ghk-bcn-90,gkk-dar-95,hl-aolnc-93};
see Kelly~\cite{k-ln-91} for an extensive survey.

Dynamic routing can be viewed as a special case of the 
on-line load balancing problem, where the load 
(incoming calls or requests) may be assigned to one 
or more servers (network links), and jobs (communication 
requests) can be scheduled only on specific subsets 
(paths) of the set of servers, as defined by the 
network topology. In this paper we study the impact 
of replacing the sequential searches of the routing 
algorithm by a version of the {\em balanced allocation 
principle}. The basic idea is as follows:
Instead of sequentially choosing alternative
options (in our case, paths) until a desirable 
one is found, in the balanced-allocation
regime the algorithm randomly chooses and 
examines a number of possible options,
and assigns the job at hand to the option
which appears to be the best at the 
time of the assignment. 

A number of papers have demonstrated 
the advantage of the application of the 
balanced allocation-principle~\cite{abku:94,abku-ba-00,
bflpr-batli-95,m-ptcrlb-96,m-arlbs-97} 
for standard load balancing
problems, where jobs require 
only one server and can be 
executed by any server in the system.
This research has shown that balanced 
allocations usually produce a very 
substantial improvement in performance,
at the cost of a small increase in overhead:
Since several alternatives are examined 
even when the first alternative would 
have been satisfactory, 
the complexity of the routing algorithm
is increased. But, as has been
shown before and as we also demonstrate
in the present context, examining even a very 
small number of alternative (thus increasing
overhead by a very small amount) can offer
great performance improvements.

The idea of employing the balanced allocation
principle to the problem of dynamic network 
routing as described in this paper 
was first explored in~\cite{lu-rncbptba-99}. 
In this context the goal is to reduce system 
congestion and minimize the 
blocking probability, that is, 
the probability that a call request is rejected. 
The main difficulty in applying and analyzing 
the balanced allocation principle in a network 
setting is in handling the dependencies imposed 
by the topology of the network.
The preliminary results in~\cite{lu-rncbptba-99} 
show that the advantage of balanced 
allocations is so significant 
that it holds even in the presence of 
a set of dependencies.

The performance of a routing protocol 
can be analyzed in a static (finite, discrete time) 
or in a dynamic (infinite, continuous time) setting.
The static case has been extensively studied 
in~\cite{lmu-olrrcn-02}, extending and 
strengthening the results in~\cite{lu-rncbptba-99}.  
In this paper we consider the continuous-time case.
The analysis of the continuous-time case 
suggested in~\cite{lu-rncbptba-99} was based on
applying Kurtz's density-dependent jump Markov 
chain technique, following the supermarket model 
analysis in~\cite{m-ptcrlb-96,m-arlbs-97}. However,
since the argument there is incomplete~\cite{l-pats-00},
we present here a different analysis. Our 
results concern the long-term behavior of
large networks employing a routing protocol
based on the balanced allocations principle.
The main tools we employ are a 
Lyapunov drift criterion used
to establish the existence of a stationary 
distribution for the BDAR routing protocol, 
and a continuous-time extension of the technique 
in~\cite{abku-ba-00}, used to analyze the 
stationary behavior of a network.

Balanced allocations have also been studied 
in the context of \emph{queueing} networks, 
where analogous results 
(under different asymptotic regimes than 
the ones in this paper) 
are obtained 
in~\cite{m-ptcrlb-96,VDK:96,martin-suhov:99,suhov-vvedenskaya:02},
among others.

\subsection{Model Description and Main Results}
In the types of networks considered in this paper,
a logical link or ``bandwidth'' is reserved
between each pair of stations, and an
alternative route is only used when this
logical link has already been exhausted.
We model such a network as the complete graph
$G=(V,E)$
with $|V|=n$ vertices (stations) and 
$|E|=N=\binom n2$ edges (links). 

The input to the system is a sequence of call 
requests, which are assumed to arrive at Poisson
times: New calls onto each link (i.e., between 
each pair of nodes) arrive according to a Poisson 
process with rate $\lambda$, all arrival streams 
being independent. Similarly, the duration of
a call is independent of all arrival times 
all other call durations,
and it is exponentially distributed with mean $1/\mu$.

The routing algorithm has to process 
the calls on-line, that is, the $t$-th request 
is either assigned a path or rejected before the
algorithm receives the ($t+1$)-th request. 
Once a call is assigned to a path, that path 
cannot be changed throughout the duration of 
the call. We assume that each edge has a capacity
of $2B$ calls, where half of this capacity is
reserved for direct links (namely it will only
be used for call requests between these two
nodes), and the other half is reserved for being
used as part of an alternative route between
two stations.

As in most of our results we consider large
networks with a number $n$ of nodes growing
to infinity, we will also assume that the 
capacity parameter $B$ may vary with $n$.
Specifically, we assume that $B=B_n$ is
nondecreasing in $n$, and we also allow 
the possibility $B=\infty$.

The goal in designing an efficient
routing protocol is to assign routes 
to the maximum possible number of call 
requests without violating the capacity 
constraints on the edges. We will compare
the performance of the following two 
protocols:

\medskip

The \emph{$d$-Dynamic Alternative Routing (DAR)
algorithm} works as follows. When a new call request 
arrives, it tries to route the call through the direct 
(one-link) path. If there is no available bandwidth 
on the direct path, then the algorithm sequentially
chooses alternative routes of length two and assigns 
the call to the first available path.
Up to $d$ such choices are made, and they are made
at random. If no possible path is found, then the 
request is rejected.

The \emph{$d$-Balanced Dynamic Alternative Routing (BDAR)
algorithm} also assigns a new call request to
the direct path if there is available bandwidth.
If not, then the algorithm chooses $d$ length-two 
alternative paths at random, and compares the maximum
load among them (where the load 
of such a path is taken to be the maximum load of
the two links on that path). Then the call is assigned
to the path with the minimum load. As before,
if there is no path with free bandwidth among 
these $d$ choices, then the call is rejected.

\medskip

The model described so far, together with one of the two
protocols above, induces a continuous-time stochastic
process describing the behavior of the network. As we 
show below, this system (for fixed $n$) converges 
to a stationary regime exponentially fast. For our 
purposes, the main performance measure is the 
minimum required bandwidth that ensures that, 
under the stationary distribution of the network, 
the blocking probability (i.e., the probability
that a new call is rejected) is appropriately
small.


In this paper our main goal is to compare the performance
of the DAR algorithm with that of BDAR. It is clear that
BDAR's performance is dominated by its performance on 
alternative (length-two) routes. Therefore,
in order to simplify the analysis, we consider
a variant of BDAR, called BDAR*, which ignores 
the direct links and services each call only 
via an alternative route, making use only of 
the~$B$ alternative connections of each edge. 
In other words, we assume that each edge has
capacity~$B$ and all of it is dedicated to 
alternative routes. We show that even though 
the BDAR* policy ignores the direct links, 
it has superior performance compared to DAR.


The following result illustrates this superiority
by exhibiting explicit asymptotic bounds on 
their bandwidth requirements.
It follows from the results in
Theorems~\ref{thm:bdar_bounded} and~\ref{thm:dar_bounded}.


\begin{theorem}
Assume that all the edges have a capacity of $2B$ links.

Under the DAR policy, edge capacity
\[B=\Omega\left(\sqrt{\frac{\ln n}{d\ln\ln n}}\right),\hspace{2.7cm}\text{as
  $n\to\infty$}\]
is necessary to ensure that a new call is not lost with high probability. 

On the other hand if we perform the BDAR* policy 
(thus ignoring the~$B$ direct links), edge 
capacity
\[B=\frac{\ln\ln n}{\ln d}+o\left(\frac{\ln\ln n}{\ln d}\right),
\hspace{3cm}\text{as
$n\to\infty$}\]
suffices to ensure that a new call is not lost with high probability.
\end{theorem}

In the above result and throughout the paper,
we say that a limiting statement holds
``with high probability'' (abbreviated ``whp.'')
if it holds with probability that is at least 
$1-1/n^c$ for some constant $c>0$. For example,
when we say that a random variable ``$X_n=O(\ln n)$
whp.'' we mean that there are positive 
constants $C$ and $c$
such that $\PR(X_n\leq C\ln n)\geq 1-1/n^c$
for all $n$ large enough. Similarly, $X_n=o(\ln n)$
whp. means that there is a $c>0$ such that,
for all $\epsilon>0$, $\PR(X_n\leq \epsilon \ln n)\geq 1-1/n^c$
for all $n$ large enough.

Note that the result of Theorem~1 is exactly analogous 
to that obtained in~\cite{lmu-olrrcn-02} in the
discrete-time case.

\section{Analysis of Balanced-Allocation Routing}

This section presents the main contribution of this paper, 
a steady state analysis of the performance of the 
BDAR* routing algorithm. The network is a complete graph
with $n$ nodes and $N=\binom n2$ undirected edges.
New calls arrive at Poisson times with rate $\lambda$
and their durations are exponentially distributed with
mean $1/\mu$, as described earlier. 
As it turns out, an important parameter in the
analysis of the network load is the ratio $\rho=\lambda/\mu$.

\subsection{Unbounded capacities}
We first analyze the maximum load on edges when the 
algorithm is used on a network with unbounded edge 
capacity, corresponding to $B=B_n=\infty$.
This model induces a continuous time Markov process
$\bfPhi=\{\Phi(t)\;:\;t\geq 0\}$, where 
$\Phi(t)=(l_1(t),l_2(t),\ldots,l_N(t))$, and
each $l_i(t)$ denotes the load, at time $t$,
of the $i$th link in the network.
As we show next, this Markov process has a stationary
distribution $\pi_n$ to which it converges exponentially
fast, regardless of the initial state of the network.
We then prove a high probability 
bound on the maximum load on any edge 
in the system under this stationary distribution.

Since we are only interested in the load of the 
alternative paths on the edges, each state of
this Markov process corresponds to the load 
on edges from a collection of length-two paths. 
We say that a vector  $x=(l_1, l_2,\dots,l_N)$
is a {\it legal state} if it corresponds to the 
load on the $N$ edges from a collection 
(possibly empty) of length-two paths.
The natural state space $\Sigma$ for our
process $\bfPhi$ is then taken to be
\[\Sigma=\left\{x=(l_1, l_2,\dots,l_N)~|~
l_i\in\mathbb{N},\ \ \mbox{$x$ is a legal state}
\right\}.\]
The process $\bfPhi$ evolves on $\Sigma$ according to the
model described above. This evolution is formalized by
the transition semigroup $\{P^t: t\ge0\}$ of $\bfPhi$,
where $P^t(x,y)$ is simply the probability that 
$\bfPhi$ is in state $y$ at time $t$ given that it
was in state $x$ at time zero,
$P^t(x,y)=\PR\{\Phi(t)=y\,|\,\Phi(0)=x\}$.

Our first result shows that $\bfPhi$ has a 
stationary (or invariant) distribution
to which it converges exponentially
fast. It is stated in terms of the ``Lyapunov
function'' $V(x)$ which is defined as 1+(total 
number of active calls in state $x$):
\begin{eqnarray}
V(x)=V(l_1,l_2,\ldots,l_N)=1+\frac{1}{2}\sum_{i=1}^N l_i.
\label{eq:V}
\end{eqnarray}

\begin{theorem}
\label{thm:unbounded}
Assume that the BDAR* algorithm is used on a
network with $n$ nodes, each of which has
infinite capacity. Then the induced Markov
process $\bfPhi$ has an invariant distribution 
$\pi_n$, and, moreover, for any initial
state $x\in \Sigma$, the distribution of $\Phi(t)$
converges to $\pi_n$ exponentially fast, namely
there is a constant $\gamma<1$, such that
$$
\sup_y |P^t(x,y)-\pi_n(y)|\leq V(x)\gamma^t,
\;\;\mbox{for all $t\geq 0$ and all $x\in\Sigma$.}
$$
\end{theorem}

\begin{proof}
Our proof uses the Lyapunov drift criterion 
for the exponential ergodicity of a continuous 
time Markov processes~\cite{st-smpIII-93,dmt-euemp-95,mt-sfltgssmp-93}.
To state our main tool we recall a few definitions, adapted
to our case of countable state space.

The {\it generator} $\mathcal{A}$ of the process $\boldsymbol{\Phi}$ 
is a linear operator on functions $F:\Sigma\to{\mathbb R}$
defined by
\[\mathcal{A}F(x)=
\lim_{h\downarrow0}\frac{\E(F(\Phi(h))\,|\,x_0=x)-F(x)}h\]
whenever the above limit exists for all $x\in\Sigma.$
The {\it explosion time} of $\boldsymbol{\Phi}$ is defined as
$$\zeta=\sup_n
J_n,$$
where \[J_0=0,\qquad J_{n+1}=\inf\{t\ge J_n: \Phi_t\neq \Phi_{J_n}\}\]
($J_0, J_1,\dots$ are the jump times of the Markov process).
We say $\bfPhi$ is \emph{nonexplosive} if
$\Pr(\zeta=\infty\,|\,\Phi_0=x)=1$ for any starting state~$x$. 


The following theorem follows from the more general results
in~\cite{mt-sfltgssmp-93,dmt-euemp-95}, specialized to the 
case of a continuous-time Markov process with a countable 
state space.

\begin{theorem}
{\em \cite{mt-sfltgssmp-93,dmt-euemp-95} }
Suppose a Markov process
evolving on a countable
state space that is nonexplosive,
irreducible (with respect to 
the counting measure on $\Sigma$)
and aperiodic.
If there exists a finite set $C\subset\Sigma$,
constants $b<\infty$, $\beta>0$ and a 
function $V:\Sigma\to[1,\infty)$, 
such that,
\begin{equation}
\label{eq:drift}
\mathcal{A}V(x)\le-\beta
V(x)+b\mathbb{I}_C(x)\qquad{x\in\Sigma}\,,
\end{equation}
then the process is positive recurrent with 
some invariant probability measure $\pi$, 
and there exist constants $\gamma<1$, $D<\infty$ 
such that 
\[\sup_y\lvert P^t(x,y)-\pi(y)\rvert\le D\, V(x)\gamma^t,
\;\;\mbox{for all $t\geq 0$ and all $x\in\Sigma$.}
\]
\end{theorem}

It is easy to verify that the process is 
$\psi$-irreducible and aperiodic,
with the maximal aperiodicity measure
$\psi$ being the counting measure on 
$\Sigma$.\footnote{This follows 
along the lines of the arguments 
in Chapters~4 and~5 of 
\cite{mt-mcss-93}. In particular,
note that all sets $\{y\}\in\Sigma$ are 
$\nu_1$-small and $P^1(x,y)>0$ for all $x,y\in\Sigma$
so that in fact $\bfPhi$ is irreducible and
strongly aperiodic.}
Also the process is nonexplosive since the 
number of new calls in a given interval has 
a Poisson distribution with a finite mean,
therefore the probability of infinite number 
of transition in a finite interval is 0.

To show that the drift criterion
(\ref{eq:drift}) can satisfied, we use
the Lyapunov function 
$V(x)$=1+(total number of active calls in state $x$)
defined in (\ref{eq:V}) above.

In order to compute $\mathcal{A}V$ we notice that when a
new call enters the system, it increases the loads of two edges by 1,
hence the value of $V$ by 1, and when a call terminates the value of $V$
decreases by 1. Therefore, new calls are generated 
with rate $\lambda N$ and calls are terminated at 
a rate $\mu (V(x)-1)$.
The probability that in a time interval~$h$ there are
2 or more new calls or terminations of calls is~$o(h)$.%
\footnote{Here and in the next expression with the notation $o(h)$ we mean
that $f$ is $o(h)$ if $\lim_{h\to0}\frac{f(h)}h=0$. In the rest of the text
$o(n)$ has the usual meaning.}
Using these observations we can compute $\mathcal{A}V$:
\begin{eqnarray*}
\mathcal{A}V(x)
&=&\lim_{h\downarrow0}\frac{V(x)+\lambda N\cdot
   h- \mu \cdot (V(x)-1)\cdot h+o(h)-V(x)}h\\
&=&\lambda N-\mu V(x)+\mu 
\end{eqnarray*}

To analyze the drift condition we distinguish between the following two cases:
\begin{itemize}
\item $x\in C$:
\[\mathcal{A}V(x)=\lambda N-\mu V(x)+\mu <-\frac{\mu V(x)}{2}+\lambda
N+\mu\]
\item $x\in C^c$ ($x$ is in the complement of $C$):
\[\mathcal{A}V(x)=\lambda N -\mu V(x)+\mu \le
\frac{\mu V(x)}{2}-\mu V(x)=
-\frac{\mu V(x)}{2}.\]
\end{itemize}
Thus, the drift condition holds
for $\beta=\mu /2$ and $b=\lambda N+\mu$.
\end{proof}

Having shown the existence of an invariant limiting distribution $\pi_n$, we
now analyze the maximum load on the edges under this distribution.

\begin{theorem}
Consider a network with~$n$ nodes, 
and let $\pi_n$ be the invariant distribution 
of the induced Markov process under the BDAR*
policy with unbounded edge capacity.
Under $\pi_n$,
the maximum number of calls in any edge 
is bounded whp. by
\[\frac{\ln\ln n}{\ln d}+o\left(\frac{\ln\ln n}{\ln
    d}\right),\qquad\text{as $n\to\infty$}.\]
\end{theorem}

\begin{proof}
In order to compute the maximum edge load under the stationary
distribution, we start observing the system at some time point and study
its transient behavior; we then use the results to deduce the properties
of the invariant distribution. In particular,
we show that there exists a constant
$T=O\left(n\frac{\ln\ln n}{\ln d}\right)$, such that for any state of
the system at time $\tau-T$
that has sufficiently large probability, whp. at time $\tau$ the maximum number 
of calls on
any edge is
\[\frac{\ln\ln n}{\ln d}+o\left(\frac{\ln\ln n}{\ln d}\right).\]

The high level idea is the following: We partition the time $T$ into
$\frac{\ln\ln n}{\ln d}+o\left(\frac{\ln\ln n}{\ln d}\right)$
periods of length $O(n)$. Roughly, we argue
that at the end of the $i$-th period, whp., for each node,
the number of incident edges with load greater than~$i$ is at
most~$\alpha_i$. The $\alpha_i$ decrease
doubly exponentially, so at the end of the last period we will be able
to deduce that there are no edges with load more than
$\frac{\ln\ln n}{\ln d}$ whp. The challenge is to handle the
dependencies, as the number of calls during some period depends on the
number of calls of the previous periods. We now proceed with the
details.

Suppose that a call routed at time $t$
is assigned to edges
$e_1$ and $e_2$.
The \emph{height} of that call at edge $e_1$ is 1 plus
$l_{e_1}(t-)$.
We define the following random variables:
\begin{itemize}
\item $L_{\ge i}^v(t)$: Number of edges incident to node $v$ with load at
least $i$ at time $t$.
\item $M_{\ge i}^v(t)$: Number of calls at edges incident to $v$ with height
greater or equal to $i$ at time $t$.
\end{itemize}

Trivially we have $L_{\ge i}^v(t)\le M_{\ge i}^v(t)$.

We define the sequence of values $\{a_i\}$ which decreases doubly
exponentially:
\begin{align*}
\alpha_{\kappa}&=\frac{(n-1)\rho}{\kappa}&\text{where
  $\kappa=e\cdot\!\sqrt[d-1]{2\rho\cdot 4^d}$}\\
\alpha_i&=\frac{2\rho\cdot 4^d\cdot a_{i-1}^d}{(n-1)^{d-1}}&\text{for $i>0$ and
   $\alpha_{i-1}\ge\sqrt[d]{\frac1{\rho}n^{d-1}\ln n}$,}\\
\alpha_{i^*}&=25\ln n&\text{$i^*$ is the smallest $i$ for which
   $a_{i-1}<\sqrt[d]{\frac1{\rho}n^{d-1}\ln n}$}\\
\alpha_{i^*+1}&=10
\end{align*}

Solving the recurrence we get for $\kappa\le i<i^*$,
\begin{equation}
\label{eq:recurrence}
\begin{split}
\alpha_{i+\kappa}&=\frac{(2\rho\cdot 4^d)^{\frac{d^i-1}{d-1}}}{\kappa^{d^i}}(n-1)
=\frac1{\sqrt[d-1]{2\rho\cdot4^d}}
\cdot\left[\frac{\sqrt[d-1]{2\rho\cdot4^d}}{\kappa}\right]^{d^i}(n-1)\\
&=\frac1{\sqrt[d-1]{2\rho\cdot4^d}}\cdot\frac{n-1}{e^{d^i}}
\end{split}
\end{equation}
and for the $i^*$
\[\alpha_{i^*-1}<\sqrt[d]{\frac1{\rho}n^{d-1}\ln n}\]
which gives
\[i^*=\frac{\ln\ln n}{\ln d}+o\left(\frac{\ln\ln n}{\ln d}\right).\]

Next we define
$T=n(i^*+2)=O\left(n\frac{\ln\ln n}{\ln d}\right)$
and an increasing sequence of points in time: let 
$t_{\kappa}=\tau-T$ and
for $i>\kappa$, $t_i=t_{i-1}+n$, so that the end of the last period,
$t_{i^*+2}$, is the current time~$\tau$.

Let $E$ denote the event ``at time $t_\kappa$ there are at most
$(1+\epsilon)N\rho$ calls in the system,''
and let $$C_i=\{\forall v\in V, t\in[t_i,T]: M_{\ge i}^v(t)\le2a_i\}.$$

We will show by induction that for $i=\kappa,\dots,i^*+1$
\[\Pr(\neg C_i\,|\,E)\le\frac{2i}{n^2}\]

For the base case ($i=\kappa$), conditioning on $E$, the expected number
of calls for
a particular node~$v$ is $(1+\epsilon)(n-1)\rho$, since each existing call 
has probability
$2/n$ to have~$v$ as an endpoint. Hence, by using the Chernoff bound
\[\Pr(\text{node~$v$ has more than $(1+\delta)(n-1)\rho$ 
calls}\,|\,E)=o(1/n^c)\]
where $\epsilon<\delta<1$ and~$c$ can be any positive constant. Therefore
\[\Pr(\neg C_\kappa\,|\,E)\le n\Pr\left(M_{\ge
     \kappa}^v>\frac{2(n-1)\rho}{\kappa}\,\bigg|\,E\right)<
\frac{2\kappa}{n^2}.\]

For the induction step we assume that
\[\Pr(\neg C_{i-1}\,|\,E)\le\frac{2(i-1)}{n^2}\]

Let~$G$ denote the event ``a new call is generated with~$v$ as an 
endpoint,'' and call~$u$ the other endpoint and~$w$ the intermediate
node of the alternative path. We have
\[\begin{split}
\label{eq:qDef}
&\Pr(\text{a new call increases $M_{\ge i}^v$}\,|\,G,C_{i-1},E)\\
&\le\Pr(\text{height of new call is $\ge i$ in either $(v,w)$ or
$(w, u)$}\,|\,G,C_{i-1},E)\\
&\parbox{\textwidth}{
\begin{flalign}
&\le\left(\frac{L_{\ge i-1}^v+L_{\ge i-1}^u}{n-1}\right)^d\notag\\
&\le\left(\frac{M_{\ge i-1}^v+M_{\ge i-1}^u}{n-1}\right)^d&
&\text{since $L_{\ge i}^v(t)\le M_{\ge i}^v(t)$}\notag\\
&\le\left(\frac{2\cdot2\alpha_{i-1}}{n-1}\right)^d\Def q_i&
&\text{from the induction hypothesis}\hspace{.5in}
\end{flalign}}
\end{split}\]
Notice that for $i=\kappa,\dots,i^*$ we have
\begin{equation}
\label{eq:boundqi}
q_i\le\frac{\alpha_i}{2\rho(n-1)}.
\end{equation}


We now define
\[F_i=\{\forall v\in V: M_{\ge i}^v(t_i)<\alpha_i\}\]
and prove Lemmas~\ref{lem:go-down} and~\ref{lem:stay-down}, that
allow us to conclude that $\Pr(\neg
C_i\,|\,E)\le\dfrac{2i}{n^2}$.
\begin{lemma}
\label{lem:go-down}
Under the inductive hypothesis
\[\Pr(\neg F_i\,|\,C_{i-1},E)\le\frac1{n^2}\]
\end{lemma}

\begin{proof}
Consider the time interval $[t_{i-1},t_i]$ and recall that
$t_i-t_{i-1}=n$.

First notice that since the duration of each call follows an exponential
distribution with parameter $\mu$, the probability that a call that
is already in the system at time~$t_{i-1}$ will remain until the end of the
interval~$t_i$ is $e^{-n\mu}$.
Hence all these calls will end before
the end of the interval with exponentially high probability. To analyze
the number of the remaining calls that
were created during the period we make use of lemma~\ref{lem:newCalls}
which completes the proof of the lemma.
\end{proof}

\begin{lemma}
\label{lem:newCalls}
Consider a period of length $\Delta$ and a given node~$v$. Conditioning on
$C_{i-1}$, the number of new calls that increased $M_{\ge i}^v$ when
they were generated, and
remained until the end of the period is less than $\alpha_i$, with
probability at least $1-\frac1{n^4}$.
\end{lemma}

\begin{proof}
Each node has $n-1$ incident links in each of which new calls are
generated with rate~$\lambda$. Conditioning on having a new request
on~$v$, $M_{\ge i}^v$ is increased with probability at most~$q_i$.
Therefore the number of calls at time~$t_i$ is
stochastically dominated by that formed by a process that generates new
calls with rate~$\lambda(n-1)q_i$ which have a duration exponentially
distributed with parameter~$\mu$. This process is the same as the
infinite server Poisson queue~\cite[page 18]{r-apmoa-92} in which the
number of calls at the end of the period is distributed according to a
Poisson distribution with rate
\[\lambda(n-1)q_i\Delta p\]
where
\[p=\int_0^{\Delta}\frac{e^{-\mu(\Delta-x)}}{\Delta}\mathrm{d}x=
\frac 1{\mu\Delta}\left(1-e^{-\mu\Delta}\right)\le\frac 1{\mu\Delta}\]
So the rate is at most $\lambda q_i(n-1)/\mu=\rho q_i(n-1)$.

We now distinguish the following two cases:
\begin{itemize}
\item[Case 1:] 
For $i\le i^*$, by using Equation~\ref{eq:boundqi} we get that the
expected number of calls at the
end of the period is at most $\alpha_i/2$
and by applying a Chernoff bound\footnote{see for
example~\cite[page 416]{r-fcp-98}.} for the Poisson
distribution, we get that the probability that the number of calls
is higher than $a_i$ is bounded by
\[\frac{e^{-\frac{\alpha_i}2}(e\frac{\alpha_i}2)^{\alpha_i}}
{\alpha_i^{\alpha_i}}=
e^{-\left(\ln2-\frac12\right)\alpha_i}\]
For $i<i^*$ we have from the definition of $\alpha_i$
\[e^{-\left(\ln2-\frac12\right)\alpha_i}=
e^{-\left(\ln2-\frac12\right)\frac{2\rho\cdot4^d\alpha_{i-1}^d}{(n-1)^{d-1}}}=
e^{-\left(\ln2-\frac12\right)\frac{2\rho\cdot4^d\frac1{\rho}n^{d-1}\ln
n}{(n-1)^{d-1}}}=o\left(\frac1{n^4}\right),
\]
while for $i=i^*$ we get
\[e^{-\left(\ln2-\frac12\right)\alpha_i}=
e^{-\left(\ln2-\frac12\right)25\ln n}=o\left(\frac1{n^4}\right).\]

\item[Case 2:] 
For $i=i^*+1$, using Equation~\ref{eq:qDef} we get that the expected
number of calls at the end of the period is at most
\[\frac{4^d\cdot\alpha_{i-1}^d}{(n-1)^d}\rho(n-1)=
\frac{(4\cdot25\ln n)^d}{(n-1)^{d-1}}\rho\]
and we get the high probability result with the Chernoff bound.
\end{itemize}
\end{proof}

\begin{lemma}
\label{lem:stay-down}
Under the inductive hypothesis
\[\Pr(\neg C_i\,|\,C_{i-1},F_i,E)\le\frac1{n^2}\]
\end{lemma}

\begin{proof}
We have:
\[\begin{split}
\Pr(\neg C_i\,|\,F_i,C_{i-1},E)
&=\frac{\Pr(\neg C_i\wedge F_i\,|\,C_{i-1},E)}{\Pr(F_i\,|\,C_{i-1},E)}\\
&\le\frac n{\Pr(F_i\,|\,C_{i-1},E)}\Pr(\exists v\in V, t_a,t_b\in[t_i,T]:\\
&\qquad\qquad M_{\ge i}^v(t_a)=\alpha_i,M_{\ge 
i}^v(t_b)=2\alpha_i,M_{\ge i}^v(t)\ge
   \alpha_i\,\forall t\in[t_a,t_b]\,\big|\,C_{i-1},E)\\
&\le\frac 
n{\Pr(F_i\,|\,C_{i-1},E)}\int_{t_a=t_i}^T\int_{t_b=t_a}^T\Pr(M_{\ge 
i}^v(t_a)=\alpha_i,\\
&\qquad\qquad M_{\ge i}^v(t_b)=2\alpha_i,M_{\ge i}^v(t)\ge
   \alpha_i\,\forall t\in[t_a,t_b]\,\big|\,C_{i-1},E)\ \mathrm{d}t_b\ 
\mathrm{d}t_a\\
&\le\frac 
n{\Pr(F_i\,|\,C_{i-1},E)}\int_{t_a=t_i}^T\int_{t_b=t_a}^T\Pr(M_{\ge 
i}^v(t_a)=\alpha_i,\\
   &\qquad\qquad M_{\ge i}^v(t_b)=2\alpha_i\,\big|\,M_{\ge i}^v(t)\ge
   \alpha_i\,\forall t\in[t_a,t_b],C_{i-1},E)\ \mathrm{d}t_b\ \mathrm{d}t_a
\end{split}\]
The probability inside the integrals is the probability that the new
calls generated during the interval $[t_a,t_b]$, increased $M_{\ge i}^v$,
and remained until the
end of the interval, are at least $\alpha_i$. By applying
Lemma~\ref{lem:newCalls}, we get that this probability is at most
$n^{-4}$. Hence
\[\begin{split}
\Pr(\neg C_i\,|\,F_i,C_{i-1},E)
&\le\frac n{1-\frac1{n^2}}\int_{t_a=t_i}^T\int_{t_b=t_a}^T\frac1{n^4}
\ \mathrm{d}t_b\ \mathrm{d}t_a\\
&\le\frac n{1-\frac1{n^2}}\cdot T^2\cdot\frac1{n^4}\\
&=o\left(\frac1{n^2}\right)
\end{split}\]
since $T=O\left(n\frac{\ln\ln n}{\ln d}\right)$.
\end{proof}

Having proven the two lemmas we can now show that $\Pr(\neg 
C_i\,|\,E)\le2i/n^2$:
\[\begin{split}
\Pr(\neg C_i\,|\,E)	&=\Pr(\neg C_i\,|\,C_{i-1},E)\cdot\Pr(C_{i-1},E)\\
			&+\Pr(\neg C_i\,|\,\neg C_{i-1},E)
			  \cdot\Pr(\neg C_{i-1},E)\\
			&\le\Pr(\neg C_i\,|\,C_{i-1},E)+\frac{2(i-1)}{n^2}\\
			&=\Pr(\neg C_i\,|\,C_{i-1},F_i,E)
			  \cdot\Pr(F_i\,|\,C_{i-1},E)\\
			&+\Pr(\neg C_i\,|\,C_{i-1},\neg F_i,E)\cdot
			  \Pr(\neg F_i\,|\,C_{i-1},E)+\frac{2(i-1)}{n^2}\\
			&\le\frac1{n^2}+\frac1{n^2}+\frac{2(i-1)}{n^2}\\
			&=\frac{2i}{n^2}
\end{split}\]

We have therefore shown that the event $C_{i^*+1}$ holds whp. until the
end of $T$, which means that for every node $v$, after the ($i^*+1$)-th
period, there will be no more than $\alpha_{i^*+1}=10$ incident edges
with load more than
$i^*+1$. We will now bound the probability that
in the next interval ($[t_{i^*+2},t_{i^*+3}]$, the last interval of $T$)
there will be an incident
edge of $v$ with load more than $i^*+3$, conditioning on the event
$C_{i^*+1}$. For this to happen, we must have at least 2 new calls to
be routed using one of the 10 high-loaded edges. The probability that
two specific new calls use these edges is at most
\begin{equation}
\label{eq:lastCalls}
\left(\frac{10}{n-2}\right)^{2d}=O\left(\frac1{n^4}\right),
\end{equation}
since $d\ge2$. The expected number of calls with $v$ as an endpoint is
$\lambda(n-1)n$, since $(n-1)$ links are connected to $v$ in each of
which new calls are generated with rate~$\lambda$, while the total length
of the interval is~$n$. This implies that whp. there will be $O(n^2)$ new
calls in the whole period. Combining this fact with
Equation~\ref{eq:lastCalls} and
summing for all the nodes we conclude that at the end of period~$T$
there will be no edges with load more than $i^*+3$ whp.

We now consider the stationary distribution $\pi_n$, and show that under
it
\[\Pr\left[l_{\max}\le\frac{\ln\ln n}{\ln d}+o\left(\frac{\ln\ln
     n}{\ln d}\right)\right]=1-o(1).\]
where $l_{\max}$ denotes the maximum number of calls on any edge.
Let $s(t)$ be the state of the system at time $\tau$, and consider the
following partitioning of the state space of the underlying Markov
process:
\begin{itemize}
\item $S_1$: States in which the total number of calls in the system is
at most $(1+\epsilon)N\rho$, and the maximum load is at most $\frac{\ln\ln
n}{\ln d}+o\left(\frac{\ln\ln n}{\ln d}\right)$.
\item $S_2$: States in which the total number of calls in the system is
at most $(1+\epsilon)N\rho$, and the maximum load is at least
$\frac{\ln\ln n}{\ln d}+\Omega\left(\frac{\ln\ln n}{\ln d}\right)$.
\item $S_3$: States in which the total number of calls in the system is
more than $(1+\epsilon)N\rho$.
\end{itemize}
We have shown that
\[\Pr(s(\tau)\in S_2\,|\,s(\tau-T)\in S_1\cup S_2)=o(1)\]
and we can easily show that
\[\Pr(s(\tau)\in S_3\,|\,s(\tau-T)\in S_1\cup S_2)=o(1)\]
Moreover in the stationary distribution the number of calls in the
system has a Poisson distribution with parameter $N$. Hence by using the
Chernoff bound
\[\sum_{i\in S_3}\pi_i=o(1)\]
Then we have
\[\sum_{i\in S_2\cup S_3}\pi_i=\sum_{i\in S_2}\pi_i+\sum_{i\in
S_3}\pi_i\]
The second term is $o(1)$, while for the first one
\begin{align*}
\sum_{i\in S_2}\pi_i&=\sum_j\Pr(s(\tau)\in S_2\,|\,S(\tau-T)=j)\cdot\pi_j\\
&=\sum_{j\in S_1\cup S_2}\Pr(s(\tau)\in S_2\,|\,S(\tau-T)=j)\cdot\pi_j\\
&\quad+\sum_{j\in S_3}\Pr(s(\tau)\in S_2\,|\,S(\tau-T)=j)\cdot\pi_j\\
&=\sum_{j\in S_1\cup S_2}\pi_j\cdot o(1)+o(1)=o(1)
\end{align*}
Therefore
\[\sum_{i\in S_2\cup S_3}\pi_i=o(1)\]
which implies that
\[\sum_{i\in S_1}\pi_i=1-o(1)\]
and completes the proof of the theorem.
\end{proof}

\subsection{Bounded Capacities}

In this section we use the analysis of the BDAR* algorithm for unbounded
capacities to compute the bandwidth requirement $B$~($<\infty$) that ensures
that a new call is not lost whp.

\begin{theorem}
\label{thm:bdar_bounded}
Assume that all the edges have capacity~$B$ circuits which can be a 
function of~$n$. Then if we perform the BDAR* policy, edge capacity
\[B=\frac{\ln\ln n}{\ln d}+o\left(\frac{\ln\ln n}{\ln d}\right),
\hspace{3cm}\text{as $n\to\infty$}\]
ensures that a new call is not lost whp.
\end{theorem}

\begin{proof}

The result for finite~$B$ follows from the proof of
Theorem~\ref{thm:unbounded} which concerns unbounded capacity.
Since the Markov process is finite and aperiodic there exists a
stationary distribution. Moreover, the analysis for the unbounded case
still holds for finite~$B$ as long as $B\le i^*+1$.


A new call will be rejected if all the~$d$ choices select one of the
edges with load $i^*+1=\ln\ln n/\ln d+o(\ln\ln n/\ln d)$.
With probability at least $1-n^2$, for each node,
the number of incident edges with load at least $i^*+1$ is at most
$2\alpha_{i^*+1}$. Therefore the probability for a call to be
rejected is no more than
\[\frac1{n^2}+\left(\frac{2\alpha_{i^*+1}}{n-1}+\frac{2\alpha_{i^*+1}}{n-1}\right)^d
=o\left(\frac1n\right)\]
since $\alpha_{i^*+1}=10$.
\end{proof}

\section{Lower Bound on the Performance of the DAR Algorithm}

To demonstrate the advantage of the balanced-allocation method
we prove here a lower bound on the maximum channel load when
requests are routed using the DAR algorithm. This bound shows an
exponential gap between the capacity required by the balanced-allocation 
algorithm and the capacity required by the standard DAR algorithm
for the same stream of inputs. Again we consider
a complete network on $n$ nodes and $N=\binom n 2$ edges.
Requests for connections between a given pair arrive according to
a Poisson process with rate $\lambda$, the duration of a connection
has an exponential distribution with expectation $1/\mu$. 


\begin{theorem}
\label{thm:dar_bounded}
Assume that all the edges have capacity~$2B$ circuits which can be a 
function of~$n$. Then if we perform the DAR policy, edge capacity
\[B=\Omega\left(\sqrt{\frac{\ln n}{d\ln\ln n}}\right),\hspace{2.7cm}\text{as
  $n\to\infty$}\]
is necessary to ensure that a new call is not lost whp.
\end{theorem}

\begin{proof}
Recall that the edges have capacities~$2B$,
capacity $B$ is used for direct connections, and the remaining capacity $B$
is used for alternative routes. 
We will compute a lower bound on the probability $P=P(B)$, that
a request arriving at an arbitrary time $t$ is rejected.

We consider first the probability $P_1$ that the new call is not
routed through the direct link. The process of routing calls through
the direct link is similar to serving customers in an
$M/M/B/B$ loss system (Poisson
arrival, exponential service time, $B$ servers, up to $B$ customers in the
system).
Applying Erlang's loss formula
(e.g.,~\cite{k-ln-91}), 
\begin{equation}
\label{eq:erlang}
P_1 = \frac{(\lambda/\mu)^B}{B!}\left(\sum_{i=0}^B
\frac{(\lambda/\mu)
^B}{i!}\right)^{-1}\geq
e^{-\lambda/\mu}\frac{(\lambda/\mu)^B}{B!}.
\end{equation}

We will now estimate the probability~$P_2$ that a request which was
generated
at time~$t$ on edge~$e$ and failed to use the direct link~$e$, fails
to be routed by an alternative path (i.e., all the~$d$ attempts to find
a non-saturated alternative path do not succeed). To give a lower bound
to the failure probability, we consider a modified system that up to
time~$t$ behaves differently from the real one by rejecting more calls
than the real one. Specifically, whenever the direct link is saturated
it tries \emph{only one} alternative path and if any of the edges of the
path are saturated the call is rejected.
Thus, clearly
more calls are lost in the modified system and therefore fewer calls
will exist at time~$t$. Notice though that from time~$t$ the system
behaves in the regular way according to the DAR algorithm.

In order to estimate the probability~$P_2$, we will try
to lower bound
the probability that at time~$t$ all the~$d$ alternative paths selected
as candidates to serve the request
have a saturated (as far as the bandwidth for alternative routes is
concerned) edge, which is lower bounded by the probability that
all the~$d$ edges selected have the corresponding edge~$e_i$
saturated (see Figure~\ref{fig:DAR}). For this we consider the system at
some prior time $t-\tau$, for some~$\tau$ that we will fix later. If some edge~$e_i$ 
at that time point is saturated with~$B$ calls, then the probability
that all these calls remain in~$e_i$ until time~$t$ equals
\[P_{\text{remain}}=e^{-\mu B\tau}.\]

Assume now that edges~$e_i$ and~$e_j$ are not saturated at
time~$t-\tau$, and let $e_{ij}$ be the edge that joins them. We will try
to compute the probability~$P_{\text{new}}$ that a request is generated
during~$\tau$ by $e_{ij}$, routed through the alternative
path~$e_i\!-\!e_j$, and remained until time~$t$. To simplify the
argument, we ignore any already existing calls from $e_{ij}$ routed
through that path---we are allowed to do that as these calls only
increase the usage of~$e_i$ and~$e_j$.
\begin{figure}[htbp]
\begin{center}
\input{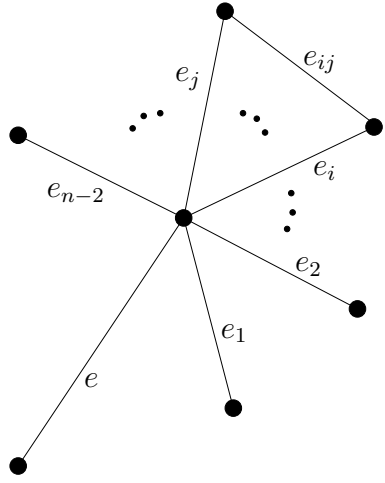}
\caption{Edge~$e$ with the new request and edges of the alternative
paths.
\label{fig:DAR}}
\end{center}
\end{figure}

We notice the following
facts:
\begin{enumerate}
\item All the direct connections of~$e_{ij}$ are occupied at time~$t-\tau$ 
with probability~$P_1$.
\item The time of a new event from edge~$e_{ij}$ (either a new call or a
termination of an existing call) is exponentially distributed with
parameter
$\lambda+B\mu$. Hence the probability of the first new event taking
place in the period~$\tau$
equals
\[1-e^{-(\lambda+B\mu)\tau}.\]
\item Conditioning that there is such a new event, the probability that
it is a new call (which will have to use an alternative path if all the
direct links are occupied) is
\[\frac{\lambda}{\lambda+B\mu}.\]
\item The probability that this call is served by the path~$e_i\!-\!e_j$
is $1/(n-2)$.
\item The probability that the call remains in the system until
time~$t$ is at least~$e^{-\mu\tau}$.
\end{enumerate}
Taking all these facts into account, we deduce that the probability that
at time~$t$ we have a call from the edge~$e_{ij}$ in~$e_i$ and~$e_j$ is
at least
\[\begin{split}
P_{\text{new}}&=P_1\cdot\left(1-e^{-(\lambda+B\mu)\tau}\right)\cdot\frac{\lambda}{\lambda+B\mu}
\cdot\frac1{n-2}\cdot e^{-\mu\tau}\\
&\ge e^{-\lambda/\mu}\frac{(\lambda/\mu)^B}{B!}
\cdot\left(1-\frac1e\right)\cdot\frac{\lambda}{\lambda+B\mu}
\cdot\frac1n\cdot e^{-\mu},
\end{split}\]
where we have selected $\tau=1$ and hence
$e^{-(\lambda+B\mu)\tau}<1/e$ for large enough~$n$.

For each edge~$e_i$ there are $n-3$ potential sources~$e_{ij}$ that are
mutually independent. Notice however that if an edge~$e_k$ is saturated
then a diverted call from edge~$e_{ik}$ that selected the
path~$e_i\!-\!e_k$ will be rejected and not contribute to the increase
of the load of~$e_i$. We perform the above counting as long as there are
at least $n/2-1$ non-saturated edges~$e_j$. Then the
probability that an edge~$e_i$ is saturated at time~$t$ is the
minimum of $P_{\text{remain}}$ and
\[P_{\text{full}}=\binom{n/2-1}{B}
P_{\text{new}}^B(1-P_{\text{new}})^{n/2-B},\]
and that minimum is always equal to $P_{\text{full}}$.
Notice that a trivial upper bound for $P_{\text{new}}$ is $1/n$ so
\[P_{\text{full}}\ge\binom{n/2-1}{B}
P_{\text{new}}^B\left(1-\frac1n\right)^{n/2-B}.\]

Let us now compute the probability~$P_2$. There are at least $n/2-1$
edges~$e_i$
whose probability of being saturated is at least~$P_{\text{full}}$,
hence the
probability that all the~$d$ alternative paths selected contain one of
the saturated edges~$e_i$ is lower bounded by
\[P_2=\left(\frac {n/2-1}{n-2}\binom{n/2-dB-1}{B}P_{\text{new}}^B
\left(1-\frac1n\right)^{n/2-B}\right)^d,\]
where the extra term $-dB$ is needed to avoid dependences between the
different edges~$e_i$. Substituting the value for $P_{\text{new}}$ we get
\[\begin{split}
P_2&\ge\frac1{2^d}\binom{n/2-dB-1}{B}^d
e^{-dB\lambda/\mu}\cdot
\frac{(\lambda/\mu)^{dB^2}}{(B!)^{dB}}
\cdot\left(1-\frac1e\right)^{dB}\cdot\\
&\quad\left(\frac{\lambda}{\lambda+B\mu}\right)^{dB}
\cdot\frac1{n^{dB}}\cdot e^{-dB\mu}\cdot\left(1-\frac1n\right)^{dn}
\end{split}\]

Therefore the probability that the call generated
at time~$t$ is rejected is at least
\[\begin{split}
P_1\cdot P_2&\ge 
e^{-\lambda/\mu}\frac{(\lambda/\mu)^B}{B!}\cdot
\frac1{2^d}\left(\frac{n/2-dB-1}{B}\right)^{dB}
e^{-dB\lambda/\mu}\cdot
\frac{(\lambda/\mu)^{dB^2}}{(B!)^{dB}}\cdot\\
&\quad\left(1-\frac1e\right)^{dB}\cdot
\left(\frac{\lambda}{\lambda+B\mu}\right)^{dB}
\cdot\frac1{n^{dB}}\cdot e^{-dB\mu}\cdot\left(\frac1{3^d}\right)\\
&=e^{-O(dB^2\ln B-dB^2\ln(\lambda/\mu))}
\end{split}\]

Therefore, in order to guarantee that
a new call is not lost whp., the bandwidth must be at
least
\[B=\Omega\left(\sqrt{\frac{\ln n}{d\ln\ln n}}\right).\]
\end{proof}

\bibliographystyle{abbrv}

\begin{thebibliography}{10}

\bibitem{acm-dondr-81}
G.~R. Ash, R.~H. Cardwell, and R.~P. Murray.
\newblock Design and optimization of networks with dynamic routing.
\newblock {\em BSTJ}, 60, 8(8):1787--1820, 1981.

\bibitem{abku:94}
Y.~Azar, A.~Broder, A.~Karlin, and E.~Upfal.
\newblock Balanced allocations.
\newblock In {\em Proceedings of the 26th ACM Symposium on the Theory of
  Computing}, pages 593--602, 1994.

\bibitem{abku-ba-00}
Y.~Azar, A.~Z. Broder, A.~R. Karlin, and E.~Upfal.
\newblock Balanced allocations.
\newblock {\em SIAM Journal on Computing}, 29(1):180--200, Feb. 2000.

\bibitem{bflpr-batli-95}
A.~Z. Broder, A.~Frieze, C.~Lund, S.~Phillips, and N.~Reingold.
\newblock Balanced allocations for tree-like inputs.
\newblock {\em Information Processing Letters}, 55(6):329--332, Sept. 1995.

\bibitem{dmt-euemp-95}
D.~Down, S.~P. Meyn, and R.~Tweedie.
\newblock Exponential and uniform ergodicity of {M}arkov processes.
\newblock {\em Ann. Probab.}, 23(4):1671--1691, 1996.

\bibitem{ghk-bcn-90}
R.~J. Gibbens, P.~J. Hunt, and F.~P. Kelly.
\newblock Bistability in communication networks.
\newblock In G.~R. Grimmet and D.~J.~A. Welsh, editors, {\em Disorder in
  Physical Systems}, pages 113--128. Oxford Univ. Press, New York, 1990.

\bibitem{gkk-dar-95}
R.~J. Gibbens, F.~P. Kelly, and P.~B. Key.
\newblock Dynamic alternative routing.
\newblock In M.~E. Steenstrup, editor, {\em Routing in Communications
  Networks}, pages 13--47. Prentice Hall, 1995.

\bibitem{hl-aolnc-93}
P.~J. Hunt and C.~N. Laws.
\newblock Asymptotically optimal loss network control.
\newblock {\em Mathematics of Operations Research}, 18(4):880--900, 1993.

\bibitem{k-ln-91}
F.~P. Kelly.
\newblock Loss networks.
\newblock {\em Annals of Applied Probability}, 1(3):319--378, 1991.

\bibitem{l-pats-00}
M.~J. Luczak.
\newblock {\em Probability, Algorithms and Telecommunication Systems}.
\newblock {DP}hil thesis, Oxford University, 2000.

\bibitem{lmu-olrrcn-02}
M.~J. Luczak, C.~McDiarmid, and E.~Upfal.
\newblock On-line routing of random calls in networks.
\newblock {\em Probability Theory and Related Fields}, 2002.
\newblock To appear.

\bibitem{lu-rncbptba-99}
M.~J. Luczak and E.~Upfal.
\newblock Reducing network congestion and blocking probability through balanced
  allocation.
\newblock In {\em {IEEE} Symposium on Foundations of Computer Science}, pages
  587--595, 1999.

\bibitem{martin-suhov:99}
J.~Martin and Y.~Suhov.
\newblock Fast {J}ackson networks.
\newblock {\em Ann. Appl. Probab.}, 9(3):854--870, 1999.

\bibitem{st-smpIII-93}
S.~P. Meyn and R.~Tweedie.
\newblock Stability of {M}arkovian processes {III}: {F}oster-{L}yapunov
  criteria for continuous-time processes.
\newblock {\em Adv. Appl. Probab.}, 25:518--548, 1993.

\bibitem{mt-sfltgssmp-93}
S.~P. Meyn and R.~Tweedie.
\newblock A survey of {F}oster-{L}yapunov techniques for general state space
  {M}arkov processes.
\newblock In {\em Proceedings of the Workshop on Stochastic Stability and
  Stochastic Stabilization, Metz, France, June 1993}. Springer-Verlag, 1994.

\bibitem{mt-mcss-93}
S.~P. Meyn and R.~L. Tweedie.
\newblock {\em {M}arkov Chains and Stochastic Stability}.
\newblock Communications and Control Engineering Series. Springer-Verlag,
  London, New York, 1993.

\bibitem{m-ptcrlb-96}
M.~Mitzenmacher.
\newblock {\em The Power of Two Choices in Randomized Load Balancing}.
\newblock PhD thesis, University of California, Berkeley, August 1996.

\bibitem{m-arlbs-97}
M.~Mitzenmacher.
\newblock On the analysis of randomized load balancing schemes.
\newblock In {\em Proceedings of the 9th Annual {ACM} Symposium on Parallel
  Algorithms and Architectures}, pages 292--301, Newport, Rhode Island, June
  22--25, 1997. SIGACT/SIGARCH and EATCS.
\newblock Extended abstract.

\bibitem{r-apmoa-92}
S.~M. Ross.
\newblock {\em Applied Probability Models with Optimization Applications}.
\newblock Dover Publications, Reprint, 1970.

\bibitem{r-fcp-98}
S.~M. Ross.
\newblock {\em A~First Course in Probability}.
\newblock Macmillan, London, 5th edition, 1998.

\bibitem{suhov-vvedenskaya:02}
Y.~Suhov and N.~Vvedenskaya.
\newblock Fast {J}ackson networks with dynamic routing.
\newblock {\em Problems of Information Transmission}, 38(2):136--153, 2002.

\bibitem{VDK:96}
N.~Vvedenskaya, R.~Dobrushin, and F.~Karpelevich.
\newblock A queueing system with a choice of the shorter of two queues -- an
  asymptotic approach.
\newblock {\em Problemy Peredachi Informatsii}, 32(1):20--34, 1996.

\end{thebibliography}

\end{document}